\documentclass{amsart}

\usepackage{amssymb}

\date{June 30, 2003}

\newcommand{\eq}[1]{~\eqref{#1}} 

\numberwithin{equation}{section}   

\newcommand{\la}{\langle}	   	
\newcommand{\ra}{\rangle}          	

\renewcommand{\(}{\left(}
\renewcommand{\)}{\right)}
\renewcommand{\[}{\left[}
\renewcommand{\]}{\right]}

\newcommand{\R}{{\mathbf R}}		 
\renewcommand{\phi}{\varphi}		
\newcommand{\cn}{\colon}	 

\catcode`@=11 \@mparswitchfalse  

\newcounter{mnotecount}[section]

\newtheorem{theorem}{Theorem}
\newtheorem*{corollary}{Corollary}

\swapnumbers

\newtheorem{thm}{Theorem}[section]
\newtheorem{lemma}[thm]{Lemma}
\newtheorem{prop}[thm]{Proposition}
\newtheorem{cor}[thm]{Corollary}
\newtheorem{claim}[thm]{Claim}
\newtheorem{mlemma}[thm]{Main Lemma}

\theoremstyle{definition}

\newtheorem*{definition}{Definition}

\theoremstyle{remark}

\newtheorem{remark}[thm]{Remark}

\newcommand{\f}{\partial}
\newcommand{\length}{\operatorname{Length}\nolimits}
\newcommand{\cd}{,\dots,}     
\newcommand{\con}{\nabla}     
\renewcommand{\setminus}{\smallsetminus}	
\newcommand{\Z}{{\mathbf Z}}
\newcommand{\C}{{\mathbf C}}
\newcommand{\e}{\varepsilon}
\newcommand{\trace}{\operatorname{tr}}
\newcommand{\ol}{\overline}
\newcommand{\s}{{\mathbb S}}
\newcommand{\interior}{\operatorname{int}}

\title{Convex Bodies of Constant Width and Constant Brightness}

\author{Ralph Howard}

\address{Department of Mathematics,
University of South Carolina,
Columbia, S.C. 29208, USA}
\email{howard\char'100math.sc.edu}
\urladdr{www.math.sc.edu/$\sim$howard}

\begin{document}


\maketitle

\section{Introduction.}

A \emph{convex body} in the $n$-dimensional Euclidean space $\R^n$ is
a compact convex set with non-empty interior.  A convex body $K$ in
three dimensional Euclidean space has \emph{constant width} $w$ iff
the orthogonal projection of $K$ onto every line is an interval of
length $w$.  It has \emph{constant brightness} $b$ iff the orthogonal
projection of $K$ onto every plane is a region of area~$b$.

\begin{theorem}\label{thm:main}
Any convex body in $\R^3$ of constant width and constant brightness is a
Euclidean ball.
\end{theorem}

Under the extra assumption that the boundary is of class $C^2$ this
was proven by S.~Nakajima (= A.~Matsumura)~\cite{Nakajima:ball} in
1926 (versions of Nakajima's proof can be found in the books of
Bonnesen and Fenchel~\cite[Sec.~68]{Bonnesen-Fenchel} and
Gardner~\cite[p.~117]{Gardner:book}).  Since then the problem of
determining if there is a non-smooth non-spherical convex body in
$\R^3$ of constant width and constant brightness has become well known
among geometers studying convexity (cf.~\cite[p.~992]{Chakerian:IQ},
\cite[Prob.~3.9 p.~119]{Gardner:book}, \cite[Ques.~2,
p.~437]{Gardner:Notices-tom}, \cite[p.~368]{Heil-Martini}).
Theorem~\ref{thm:main} solves this problem.

For convex bodies with $C^2$ boundaries and positive curvature
Nakajima's result was generalized by
Chakerian~\cite{Chakerian:rel-width} in 1967 to ``relative geometry''
where the width and brightness are measured with with respect to some
convex body $K_0$ symmetric about the origin called the \emph{gauge
body}.  While the main result of this paper is Theorem~\ref{thm:main},
Chakerian's methods generalize and simplify parts of our original
proof.  The following isolates the
properties required of the gauge body.  Recall the \emph{Minkowski
sum} of two subsets $A$ and $B$ of $\R^n$ is $A+B=\{a+b:a\in A, b\in
B\}$.

\begin{definition}\label{def:gauge}
A convex body $K_0$ is a \emph{regular gauge} iff it is centrally
symmetric about the origin and there are convex sets $K_1$, $K_2$ and
Euclidean balls $B_r$ and $B_R$ such that $K_0=K_1+B_r$ and
$B_R=K_0+K_2$. 
\end{definition}

Any convex body symmetric about the origin with $C^2$ boundary and
positive Gaussian curvature is a regular gauge
(Corollary~\ref{C2-gauge} below).  For any linear subspace $P$ of
$\R^n$ let $K|P$ be the projection of $K$ onto $P$ (all projections in
this paper are orthogonal).  For a unit vector $u$ let $w_K(u)$ be
the width in the direction of $u$.  For each positive integer $k$ and
any Borel subset $A$ of $\R^n$ let be $V_k(A)$ be the $k$-dimensional
volume of $A$ (which in this paper is the $k$-dimensional Hausdorff
measure of $A$).  Two subsets $A$ and $B$ of $\R^n$ are
\emph{homothetic} iff there is a positive scalar $\lambda$ and a
vector $v_0$ such that $B=v_0+\lambda A$.

\begin{theorem}\label{main:rel}
Let $K_0$ be a regular gauge in $\R^3$ and let $K$ be any convex body
in $\R^3$ such that for some constants $\alpha$, $\beta$ the
equalities $ w_K(u)=\alpha w_{K_0}(u)$ and $V_2(K|u^\bot)=\beta
V_2(K_0|u^\bot)$ hold for all $u\in \s^{n-1}$.  Then $K$ is homothetic
to $K_0$.
\end{theorem}

Letting $K_0$ be a Euclidean ball recovers Theorem~\ref{thm:main}.
While we are assuming some regularity on the gauge body $K_0$, the
main point is that no assumptions, other than convexity, are being put
on $K$.  It is likely that the result also holds with no restrictions
on either $K$ or $K_0$.  One indication this may be the case is a
beautiful and surprising result of Schneider~\cite{Schneider:poly}
that almost every, in the sense of Baire category, centrally convex
body $K_0$ is determined up to translation in the class of all convex
bodies by just its width function.  This contrasts strongly with the
fact that for any regular gauge $K_0$ there is an infinite dimensional
family of convex bodies that have the same width function as $K_0$
(see Remark~\ref{rmk:width} below).

Two convex bodies $K$ and $K_0$ in $\R^n$ have \emph{proportional
$k$-brightness} iff there is a constant $\gamma$ such that
$V_k(K|P)=\gamma V_k(K_0|P)$ for all $k$-dimensional subspaces $P$ of
$\R^n$.  Theorem~\ref{main:rel} implies a result, valid in all
dimensions, about pairs of convex bodies that have both $1$-brightness and
$2$-brightness proportional.  If $A$ and $B$ are convex sets in $\R^n$
and $L$ is a linear subspace of $\R^n$, then taking Minkowski sums
commutes with projection onto $L$, that is $(A+B)|L=A|L+B|L$.  As the
projection of a Euclidean ball is a Euclidean ball, it follows that if
$K_0$ is a regular gauge in $\R^n$, then $K_0|L$ is a regular gauge in
$L$.  Also, if $P$ is a linear subspace of $L$, then $K|P=(K|L)|P$.
Therefore if $K_0$ is a regular gauge in $\R^n$ and $K$ is a convex
body such that $K$ and $K_0$ have proportional $1$-brightness and
proportional $2$-brightness, then for any three dimensional subspace
$L$ of $\R^n$ the set $L|K_0$ is a regular gauge in $L$ and $K_0|L$
and $K|L$ will have proportional $1$-brightness and proportional
$2$-brightness as subsets of $L$.  Thus by Theorem~\ref{main:rel}
$K|L$ is homothetic to $K_0|L$.  However, if the projections $K_0|L$
and $K|L$ are homothetic for all three dimensional subspaces $L$,
then, \cite[Thm~3.1.3, p.~93]{Gardner:book}, $K$ is homothetic to
$K_0$.  Thus:

\begin{corollary}\label{cor:main}
If $K_0$ is a regular gauge in $\R^n$, $n\ge 3$, and $K$ is a convex
body in $\R^n$ that has $1$-brightness and $2$-brightness proportional
to those of $K_0$, then $K$ is homothetic to $K_0$.  In particular if
$K_0$ a Euclidean ball this implies any convex body $K$ in $\R^n$ of
constant $1$-brightness and $2$-brightness is also a Euclidean
ball.\qed 
\end{corollary}

The contents of this the paper are as follows.  In Section~\ref{sec:prelim}
some preliminaries about convex sets are given and a $C^{1,1}$
regularity result, Proposition~\ref{C11-sum}, for the support
functions of convex sets in $\R^n$ that appear is a summand in a
convex set with $C^{1,1}$ support function is proven. (I am indebted
to Daniel Hug for some of the results in this section).
Section~\ref{sec:spt} gives explicit formulas, in terms of the support
function, $h$, for the inverse of the Gauss map of the boundary of a
convex set in $\R^n$ and conditions are given for two convex sets with
$C^{1,1}$ boundary to have proportional brightness.  It is
important for our applications that some of these formulas
(eg.~Proposition~\ref{Gauss-def}) apply even when the function $h$ is
not the support function of a convex set.  In Section~\ref{sec:3D} the
results of the previous sections are applied to reduce the proof
Theorem~\ref{main:rel} to an analytic problem.  In
Section~\ref{sec:quasi} the analytic result is proven by use of
quasiconformal maps, the Beltrami equation, and the elementary theory
of covering spaces.

\section{Preliminaries on convexity.}\label{sec:prelim}

We assume that $\R^n$ has its standard inner product $\la\,,\ra$ and
let $\s^{n-1}$ be the unit sphere of $\R^n$.  For any convex body $K$
contained $\R^n$, the support function $h=h_K$ of $K$ is the function
$h\cn \s^{n-1}\to \R$ given by $ h(u):=\max_{y\in K}\la y,u\ra$. A
convex body is uniquely determined by its support function.  The
Minkowski sum of $K_1$ and $K_2$ corresponds to the sum of the support
functions: $h_{K_1+K_2}= h_{K_1}+h_{K_2}$.  The \emph{width function}
of $K$ is $w=w_K$ is $ w(u)=h(u)+h(-u)$.  This is the length of the
projection of $K$ onto a line parallel to the vector $u$.  In the
terminology of Gardner,~\cite[p.~99]{Gardner:book}, the \emph{central
symmetral} of a convex body $K$ is the convex body
$K_0:=\frac12(K-K)=\{\frac12(a-b): a,b\in K\}$.  The body $K_0$ is
centrally symmetric about the origin, and, denoting the support
function of $K_0$ by $h_0$, it follows from
$h_{\frac12(K-K)}=\frac12h_K+\frac12h_{-K}$ that $
h_0(u)=\frac12(h(u)+h(-u))$.  Therefore $K$ and $K_0$ have the same
width in all directions.  These definitions imply that a convex body
has constant width $w$ if and only if its central symmetral is a
Euclidean ball of radius $w/2$.

We need the following, which is an elementary corollary of the
Brunn-Minkowski theorem. For a proof see~\cite[Thm~3.2.2,
p.~100]{Gardner:book}. 

\begin{prop}\label{K0-vol}
The volumes of a convex body $K$ and its central symmetral
$K_0=\frac12(K-K)$ satisfy $V(K_0)\ge V(K)$ with equality if
and only if $K$ is a translate of $K_0$.\qed
\end{prop}

Recall that a function $f$ defined on an open subset $U$ of $\R^k$ is
of class $C^{1,1}$ iff it is continuously differentiable and all the
first partial derivatives satisfy a Lipschitz condition.  A convex
body $K$ has $C^{1,1}$ boundary iff its boundary $\f K$ is locally the
graph of a $C^{1,1}$ function.

There is a very nice geometric characterization of the convex bodies
that have $C^{1,1}$ boundaries in terms of freely sliding bodies.  Let
$K_1$ and $K_2$ be convex bodies in $\R^n$.  Then $K_1$ \emph{slides
freely inside of} $K_2$ iff for all $a\in \f K_1$ there is a translate
$y+K_2$ of $K_2$ such that $K_1\subseteq y+K_2$ and $a\in y+K_2$.  It
is not hard to see,~\cite[Thm~3.2.2, p.~143]{Schneider:convex}, that
$K_1$ slides freely inside of $K_2$ if and if $K_1$ is a Minkowski
summand of $K_2$.  That is, if and only if there is a convex set $K$
such that $K+K_1=K_2$.  In what follows we will use the expressions
``$K_1$ slides freely inside of $K_2$'' and ``$K_1$ is a Minkowski
summand of $K_2$'' interchangeably.  A proof of the following can be
found in~\cite[Prop.~1.4.3, p.~97]{Hormander:convexity}.

\begin{prop}\label{C11-bdry}
A convex body $K$ has $C^{1,1}$ boundary if and only if some Euclidean
ball $B_r$ slides freely inside of $K$.\qed
\end{prop}

I learned of the following elegant dual from of this theorem, with a
somewhat different proof, from Daniel Hug.

\begin{prop}[D. Hug~\cite{Hug:slide}]\label{Hug:slide}
The support function $h$ of a convex body $K$ is $C^{1,1}$ if and only
if $K$ slides freely inside of some  Euclidean ball $B_R$.
\end{prop}

\begin{proof}
Assume that $K$ slides freely inside of the ball $B_R$ of radius $R$.
Without loss of generality it may be assumed that the origin is in the
interior of $K$.  Let $K^\circ:=\{y: \la y,x\ra\le 1 \text{ for all }
x\in K\}$ be the \emph{polar body} of $K$.  The radial function of
$K^\circ$ (which is the positive real valued function $\rho$ on
$\s^{n-1}$ such that $u\mapsto \rho(u)u$ parameterizes the boundary
$\f (K^\circ)$ of $K^\circ$) is $\rho(u)=1/h(u)$, \cite[Rmk~1.7.7,
p.~44]{Schneider:convex}.  So it is enough to show that $\rho$ is a
$C^{1,1}$ function, and to show this it is enough to show that the
boundary $\f (K^\circ)$ is $C^{1,1}$.  By Proposition~\ref{C11-bdry}
it is enough to show that some ball slides freely inside of $K^\circ$.
Let $\rho(u)u\in \f (K^\circ)$.  Becasue $K$ slides freely inside a
ball of radius $R$ there is a ball $B(a,R)$ of radius centered at some
point $a$ such that $K\subset B(a,R)$ and a point $x\in K\cap \f
B(a,R)$ such that $u$ is the outward pointing normal to $B(a,R)$ at
$x$.  As the operation of taking polars is inclusion reversing,
$B_R(a)^\circ$ is contained in $K^\circ$ and as $u$ is the outward
pointing unit normal to both $K$ and $B(a,R)$ at $x$ we also have
$\rho(u)u\in \f (B_R(a)^\circ)$.  The support function of $B_R(a)$ is
$h_{B_R(a)}(u)=R+\la a,u\ra$ and therefore the radial function of the
polar $B_R(a)^\circ$ is $\rho_{B_R(a)^\circ}(u)=1/(R+\la a,u\ra)$.
Thus points on $\f (B_R(a)^\circ)$ are of the form $y=\(1/(R+\la
a,u\ra)\)u$ for $u\in \s^{n-1}$.  This implies $|y|=1/(R+\la a,u\ra)$
and $\la a,y\ra=\la a,u\ra/(R+\la a,u\ra)$.  If $\la a,u\ra$ is
ellimated from these equations the result can be written as
$$
R^2|y|^2-\la a,y\ra^2+2\la a,y\ra=1.
$$
For each $a$ this is an ellipsoid and an ellipsoid has positive
rolling radius (which is the largest number $r$ so that a ball of
radius $r$ slides freely inside of the body).  More generally for any
ball $B_R(v)$ of radius $R$ and center $v$ containing $K$ the polar
$B_R(v)^\circ$ is an ellipsoid.  By Blaschke's rolling theorem,
\cite[Cor.~3.2.10, p. 150]{Schneider:convex}, the rolling radius is
the smallest radius of curvature of $\f (B_R(v)^\circ)$ and this is a
continuous function of the vector $v$.  The set of $v$ such that
$B_R(v)$ contains $K$ is a compact set and therefore, by the
continuous dependence of the rolling radius of $\f (B_R(v)^\circ)$ on
$v$, there is a positive number $r_0$ such that a
ball of radius $r_0$ slides freely inside of any $B_R(v)^\circ$ that
contains $K$.  In particular this is true of $B_R(a)^\circ$ and so
$K^\circ$ contains an internally tangent ball of radius $r_0$ at
$\rho(u)u$.  But $\rho(u)u$ was an arbitrary point of $\f (K^\circ)$
and whence a ball of radius $r_0$ slides freely inside of $K^\circ$ as
required.

Conversely assume that the support function $h$ of $K$ is $C^{1,1}$.
Let $\tilde{h}$ be the extension of $h$ to $\R^n$ that is homogeneous of
degree $1$.  Explicitly
\begin{equation}\label{extend1}
\tilde{h}(x):=\max_{y\in K}\la y,x\ra.
\end{equation}
As $h$ is $C^{1,1}$ the function $\tilde{h}$ is $C^{1,1}_\text{\rm
Loc}$ on $\R^n\setminus \{0\}$ and $\tilde{h}$ is convex,
\cite[Thm~1.7.1, p.~38]{Schneider:convex}, the distributional Hessian
$\f^2\tilde{h}$ will be positive semi-definite on $\R^n\setminus
\{0\}$ and, because $h$ is $C^{1,1}$, locally bounded above.  Thus
there is a positive real number $R$ such that $H_0:=R\|\cdot
\|-\tilde{h}$ is a convex function.  But then,~\cite[Thm~1.7.1,
p.~38]{Schneider:convex}, $H_0\big|_{\s^{n-1}}$ is the support
function of a unique convex body $K_0$ and $H_0+\tilde{h}=R\|\cdot \|$
implies that $K+K_0=B_R(0)$.  Therefore $K$ is a summand in a ball.
\end{proof}

\begin{cor}\label{C2-gauge}
Let $K_0$ be a body that is centrally symmetric about the origin, with
$\f K_0$ of class $C^2$ with positive Gauss curvature.  Then $K_0$
is a regular gauge.
\end{cor}

\begin{proof}
It follows from a generalization Blaschke's rolling theorem,
\cite[Cor. 3.2.10, p. 150]{Schneider:convex}, that if $B_r$ is a
Euclidean ball with $r$ smaller than any of the radii of curvature of
$\ K_0$, that $B_r$ slides freely inside of $K_0$ and if $R$ is larger
than any of the radii of curvature of $\f K_0$, then $K_0$ slides
freely inside of $B_R$.
\end{proof}

\begin{prop}\label{C11-sum}
Let $K_1\cd K_k$ be convex bodies in $\R^n$ such that the Minkowski
sum $K_1+\dots +K_k$ has $C^{1,1}$ support function.  Then each summand
$K_j$ also has $C^{1,1}$ support function.  
\end{prop}

\begin{proof}
If $K_1+\dots +K_k$ has $C^{1,1}$ support function then, by
Proposition~\ref{Hug:slide}, $K_1+\dots+K_k$ is a Minkowski summand in
some ball $B_R$.  But then each $K_j$ is also a summand in $B_R$
and therefore Proposition~\ref{Hug:slide} yields that $K_j$ has
$C^{1,1}$ support function.
\end{proof}

\begin{cor}\label{even-reg}
Let $K$ be a convex body such its central symmetral has 
a $C^{1,1}$ support function.  Then the support function of
$K$ is also $C^{1,1}$.  In particular any convex body of constant
width has $C^{1,1}$ support function.
\end{cor}

\begin{proof}
If $K_0$ is the central symmetral of $K$, then $K+(-K)=2K_0$.  As
$K_0$ has $C^{1,1}$ support function, $h_0$, the support function,
$2h_0$, of $2K_0$ is also $C^{1,1}$ and therefore the support function
of $K$ is $C^{1,1}$ by Proposition~\ref{C11-sum}.
\end{proof}

\begin{remark}\label{rmk:width}
Corollary~\ref{even-reg} is sharp in the sense that even when the
support function, $h_0$, of the central symmetral, $K_0$, is
$C^\infty$, the most that can be said about the regularity of support
function, $h$, of $K$ is that it is $C^{1,1}$.  For example let $h_0$
be the support function of a regular gauge, $K_0$, and let $p$ a
$C^{1,1}$ function $p\cn\s^{n-1}\to \R$ with $p(-u)=-p(u)$.  Then for
sufficiently small $\e>0$ the function $h:=h_0+\e p$ is the support
function of a convex body with the same width function as $K_0$.  But
there are many choices of $h_0$ and $p$ with $h_0$ of class $C^\infty$
and $h$ only of class $C^{1,1}$.
\end{remark}

\section{Support Functions and  the Inverse of the Gauss Map.}\label{sec:spt}

We view vector fields $\xi$ on subsets of $U$ of $\R^n$ as functions
$\xi\cn U\to \R^n$.  A vector field on $\s^{n-1}$ is a function
$\xi\cn\s^{n-1}\to \R^n$ such that for all $u\in \s^{n-1}$ the vector
$\xi(u)\in T_u\s^{n-1}$.  As the tangent space, $T_u\s^{n-1}$, to
$\s^{n-1}$ at $u$ is just $u^\bot$, the orthogonal compliment to $u$
in $\R^n$, a vector field $\xi$ on $\s^{n-1}$ can also be viewed as a
map from $\s^{n-1}$ to $\R^n$ with $\xi(u)\bot u$ for all $u$.  If
$X\in T_u\s^{n-1}$ is a tangent vector to $\s^{n-1}$ at $u$, then a
\emph{curve fitting $X$} is a smooth curve $c\cn(a,b)\to \s^{n-1}$
defined on an interval about $0$ with $c(0)=u$ and $c'(0)=X$.
If $\xi$ is a vector field on $\s^{n-1}$ that is differentiable at the
point $u$, then for any $X\in T_u\s^{n-1}$ the \emph{covariant
derivative, $(\con_X\xi)(u)$, of $\xi$ by $X$} is the 
projection of $\left.\frac{d}{dt}\xi(c(t))\right|_{t=0}$ onto
$T_u\s^{n-1}$ where $c$ is any curve fitting $X$.  This is independent
of the choice of $c$ fitting $X$ and is given explicitly by
$$
(\con_X\xi)(u):=\left.\frac{d}{dt}\xi(c(t))\right|_{t=0}
-\left\la \left.\frac{d}{dt}\xi(c(t))\right|_{t=0},u\right\ra u.
$$
This definition implies that for any smooth curve
$c\cn (a,b)\to\s^{n-1}$ and any vector field $\xi$ on $\s^{n-1}$ that
\begin{equation}\label{curve-diff}
\dfrac{d}{dt}\xi(c(t))=(\con_X\xi)(c(t))+\left\la
\frac{d}{dt}\xi(c(t)), c(t)\right\ra c(t)
\end{equation}
for any value $t$ such that $\xi$ is differentiable at $c(t)$.  

For any $C^1$ function $p\cn \s^{n-1}\to \R$ the (spherical)
\emph{gradient} is the vector field, $\con p$, on $\s^{n-1}$ such that
$\la \con p, X\ra = dp(X)$ for all vectors $X$ tangent to $\s^{n-1}$.
At any point $u$ where the vector field $\con p$ is differentiable
the \emph{second derivative} of $p$ is the linear map
$\con^2p(u)\cn T_u\s^{n-1}\to T_u\s^{n-1}$ given by
$$
\con^2p(u)X:=(\con_X\con p)(u).
$$

\begin{remark}\label{rmk:sym}
There is a another way of viewing $\con^2 p$ that is useful.  If
$p$ is defined on $\s^{n-1}$ then extend $p$ to $\R^n$ to be
homogeneous of degree one.  That is let $\tilde{p}\cn \R^n\to\R$ be
\begin{equation}\label{extend2}
\tilde{p}(x)=|x|p(|x|^{-1}x)
\end{equation}
for $x\ne0$ and $\tilde{p}(0)=0$.  Let $\f \tilde{p}$ be the usual
gradient of $\tilde{p}$, that is $\f\tilde{p}$ is the column vector with
components $\f_1\tilde{p},\f_2\tilde{p}\cd \f_n\tilde{p}$, and let
$\f^2\tilde{p}$ be the field of linear maps on $\R^n\setminus\{0\}$
given by $ \f^2\tilde{p}(x)Y:=(\f_Y\f\tilde{p})(x)$ where $\f_Y$ is the
usual directional derivative in the direction of the vector $Y$.  The
matrix of $\f^2\tilde{p}$ with respect to the coordinate basis is the
usual Hessian matrix $\[\f_i\f_j\tilde{p}\]$.  A strightforward
calculation shows that $\f^2\tilde{p}$ and $\con^2p$ are related by
\begin{equation}\label{f2-con2}
\f^2\tilde{p}(x)Y=\frac1{|x|}\Big(\con^2p(|x|^{-1}x)+p(|x|^{-1}x)I\Big)(Y-|x|^{-2}\la
Y,x\ra x).
\end{equation}
This implies that if $u\in \s^{n-1}$ and $Y\in T_u\s^{n-1}=u^\bot$,
then 
$$
\f^2\tilde{p}(u)Y=(\con^2p(u)+p(u)I)Y
$$
and $\f^2\tilde{p}(u)u=0$.  Thus $T_u\s^{n-1}$ is invariant under
$\f^2\tilde{p}$.  The symmetry of the second partials implies that
when $p$ is $C^2$, so that $\tilde{p}$ is $C^2$ on
$\R^n\setminus\{0\}$, then $\f^2\tilde{p}(x)$ is self-adjoint (that is
$\la \f^2\tilde{p}(x)X,Y\ra=\la X,\f^2\tilde{p}(x)Y\ra$) for $x\in
\R^n\setminus\{0\}$.  But then
$\con^2p(u)=\f^2\tilde{p}(u)\big|_{T_u\s^{n-1}}-p(u)I$ implies that
$\con^2p(u)$ is self-adjoint on $T_u\s^{n-1}$.  The
formula\eq{f2-con2} also implies that $\con^2p$ exists at $u\in
\s^{n-1}$ if and only if $\f^2\tilde{p}$ exists at all points $tu$
with $t>0$.  This, combined with Fubini's Theorem, yields that
$\con^2p$ exists almost everywhere on $\s^{n-1}$ if and only if
$\f^2\tilde{p}$ exists almost everywhere on $\R^n$.
\end{remark}

\begin{prop}\label{Gauss-def}
Let $\phi\cn \s^{n-1}\to \R^n$ be a Lipschitz map such that for all
$u$ where the derivative $\phi'(u)$ exists it satisfies $\phi'(u)X\in
T_u\s^{n-1}$ for all $X\in T_u\s^{n-1}$.  Then there is a unique
$C^{1,1}$ function $p\cn \s^{n-1}\to \R$ such that
\begin{equation}\label{phi-def}
\phi(u)=p(u)u+\con p(u).
\end{equation}
The derivative $\phi'(u)$ exists at $u$ if and only if the second
derivative $\con^2p(u)$ of $p$ exists at $u$ and at these points
\begin{equation}\label{phi'}
\phi'(u)=p(u)I+\con^2p(u)
\end{equation}
where $I$ is the identity map on $T_u\s^{n-1}$.  Conversely if $p$ is
$C^{1,1}$ and $\phi$ is given by~\ref{phi-def} then $\phi'(u)X\in
T_u\s^{n-1}$ for all $X\in T_u\s^{n-1}$ for all points $u$ where
$\phi$ is differentiable.  Finally for $k\ge1$ the function $\phi$ is
$C^{k}$ if and only if $p$ is $C^{k+1}$.
\end{prop}

\begin{proof}
Any function $\phi\cn \s^{n-1}\to \R^n$ can be uniquely written as $
\phi(u)=p(u)u+\xi(u)$ where $p\cn \s^{n-1}\to \R$ and $\xi$ is a
vector field on $\s^{n-1}$.  Because $\phi$ is Lipschitz, so are $p$
and $\xi$. Therefore a theorem of Rademacher,~\cite[Thm~3.1.6,
p.~216]{Federer:book}, implies that $p$ and $\xi$ are both
differentiable almost everywhere on $\s^{n-1}$.  Let $E$ be the set of
points where both $p$ and $\xi$ are differentiable.  Then $\phi$ is
also differentiable at $u$.  Let $u\in E$, $X\in T_u\s^{n-1}$, and $c$
a curve fitting $X$.  Then, using\eq{curve-diff},
\begin{align*}
\phi'(u)X&=\left.\frac{d}{dt}\right|_{t=0}\big(p(c(t))c(t)+\xi(c(t))\big)\\
&=dp_u(X)u+p(u)X+(\con_X\xi)(u)
+\left\la \left.\frac{d}{dt}\right|_{t=0}\xi(c(t)),u\right\ra u.
\end{align*}
But $dp_u(X)=\la \con p(u),X\ra$ and, using that $\la
\xi(c(t)),c(t)\ra\equiv0$,
\begin{align*}
\left\la \left.\frac{d}{dt}\right|_{t=0}\xi(c(t)),u\right\ra
&=\left.\frac{d}{dt}\right|_{t=0}\la \xi(c(t)),c(t)\ra
-\la \xi(c(t)),c'(t)\ra\bigg|_{t=0}\\
&=-\la \xi(u),X\ra.
\end{align*}
Therefore the formula for $\phi'(u)X$ becomes
\begin{equation}\label{pre-phi'}
\phi'(u)X=\la \con p(u)-\xi(u),X\ra u + p(u)X+(\con_X\xi)(u).
\end{equation}
As $\phi'(u)X\in T_u\s^{n-1}$ the component normal to $\s^{n-1}$ must
vanish.  Whence $\la \con p(u)-\xi(u),X\ra=0$ for all $X\in
T_u\s^{n-1}$.  This implies
\begin{equation}\label{xi-p}
\xi(u)=\con p(u) \quad \text{at points $u$ where both $p$ and $\xi$ are
differentiable.} 
\end{equation}
We now argue that $p$ is continuously differentiable and that $\con
p=\xi$ on all of $\s^{n-1}$.  This will be based on the following
elementary lemma, whose proof will be given after the proof of
Proposition~\ref{Gauss-def}.

\begin{lemma}\label{p-diff}
Let $q$ be a real valued Lipschitz function defined on an open subset
$U$ of $\R^N$.  Assume that there are Lipschitz functions $q_1\cd q_N$
on $U$ and a set of full measure $S\subseteq U$ such that for all
$x\in S$ the partial derivatives of $q$ exist and satisfy $\f_j
q(x)=q_j(x)$ for all $x\in S$.  Then $q$ is of class $C^{1,1}$ and
$\f_jq=q_j$ on all of $U$.
\end{lemma}

Near any point, $u_0$, of $\s^{n-1}$ there is a $C^\infty$
parameterization $f\cn U\to V\subset\s^{n-1}$ of a neighborhood $V$ of
$u_0$, with $U$ a bounded open set in $\R^{n-1}$, and $f$ a $C^\infty$
diffeomorphism.  To show that $p$ is $C^{1,1}$ it is enough to show
the function $q\cn U\to \R$ given by $q(x):=p(f(x))$ is $C^{1,1}$.
Let $S$ be the subset of points $x\in U$ where both $p$ and $\xi$ are
differentiable at $f(x)$.  As $p$ and $\xi$ are Lipschitz and $f$ is a
diffeomorphism this is a set of full measure in $U$ and at all points
of $S$ we have, by\eq{xi-p}, that $\con p(f(x))=\xi(f(x))$.  As $\xi$
is Lipschitz there are real valued Lipschitz functions $\xi^1\cd
\xi^{n-1}$ defined on $U$ such that $
\xi(f(x))=\sum_{i=1}^{n-1}\xi^i(x)\f_if(x)$.  Therefore at points $x$
in $S$ we have $ \con
p(f(x))=\xi(f(x))=\sum_{i=1}^{n-1}\xi^i(x)\f_if(x)$ and thus
$$
\f_jq(x)=dp_{f(x)}(\f_jf)=\la \con p(f(x)),\f_jf\ra
=\sum_{i=1}^{n-1}\xi^i(f(x))\la \f_if(x),\f_jf(x)\ra.
$$
The functions $q_j:=\sum_{i=1}^{n-1}\xi^i(f(x))\la
\f_if(x),\f_jf(x)\ra$ are Lipschitz so Lemma~\ref{p-diff}
implies that $q$, and therefore also $p$, is a $C^{1,1}$ function and
that $\con p$ is a Lipschitz.

By\eq{xi-p} $\con p(u)=\xi(u)$ on the dense set $E$ and $\con p$ and
$\xi$ are continuous thus $\con p=\xi$ on all of $\s^{n-1}$.
Therefore $\phi(u)$ is given by\eq{phi-def} as required.  When $\phi$
is of this form it is clear that $\phi$ is differentiable exactly at
the points $u$ where the second derivative $\con^2p(u)$ exists.  At
such points use $\con p=\xi$ and $\con_X\xi(u)=(\con_X\con p)(u)=
\con^2p(u)$ in\eq{pre-phi'} to see that\eq{phi'} holds.  This
completes the proof that if $\phi\cn \s^{n-1}\to\R^n$ is a Lipschitz
map with $\phi'(u)X\in T_u\s^{n-1}$ for all $u\in \s^{n-1}$ where
$\phi$ is differentiable, then $\phi$ is given by\eq{phi-def} for a
uniquely determined $C^{1,1}$ function $p$.

Conversely if $p$ is $C^{1,1}$ let $\xi=\con p$ in the calculations
leading up to\eq{pre-phi'} to see that $\phi$ given by\eq{phi-def}
satisfies $\phi'(u)X\in T_u\s^{n-1}$ for all $u\in \s^{n-1}$ where
$\phi$ is differentiable.

Finally $\phi(u)=h(u)u+\con p(u)$ makes it clear that if $h$ is
$C^{k+1}$, then $\phi$ is $C^k$.  Conversely if $\phi$ is $C^{k}$,
then $p(u)=\la u,\phi(u)\ra$ implies $p$ is $C^k$.  Then
$\con p(u)=\phi(u)-p(u)u$ implies that $\con p$ is also $C^k$.  But if
$\con p$ is $C^k$, then $p$ is $C^{k+1}$.
\end{proof}

\begin{proof}[Proof of Lemma~\ref{p-diff}]
We will show that the $j$-th distributional derivative of $q$ is
$q_j$.  By definition this means we need to show that for all
$C^\infty$ functions $\psi$ with compact support contained in $U$ that
$\int_U q\f_j\psi\,dx=-\int_U q_j\psi\,dx$.  Let $e_j$ be the $j$-th
coordinate vector.  Then
\begin{align*}
\int_Uq(x)\f_j\psi(x)\,dx &=\lim_{h\to
0}\int_Uq(x)\frac{\psi(x+he_j)-\psi(x)}{h}\,dx\\
&=\lim_{h\to0}\(\frac1h\int_Uq(x)\psi(x+he_j)\,dx-\frac1h\int_Uq(x)\psi(x)\,dx\)\\
&=\lim_{h\to0}\(\frac1h\int_Uq(x-he_j)\psi(x)\,dx-\frac1h\int_Uq(x)\psi(x)\,dx\)\\
&=\lim_{h\to 0}\int_U\frac{q(x-he_j)-q(x)}{h}\psi(x)\,dx.
\end{align*}
But $q$ is Lipschitz and therefore the quotients $(q(x-he_j)-q(x))/h$
are uniformly bounded.  By assumption for all $x\in S$,
$\lim_{h\to0}(q(x-he_j)-q(x))/h=-\f_jq(x)=-q_j(x)$ and $S$ has full
measure so this limit holds almost everywhere.  Therefore Lebesgue's
bounded convergence theorem implies $\lim_{h\to
0}\int_U((q(x-he_j)-q(x))/{h})\psi(x))\,dx= -\int_U
q_jq(x)\psi(x)\,dx$.  Using this in the calculation above yields that$
\int_U q\f_j\psi\,dx=-\int_U q_j\psi\,dx$ holds, and thus the
distributional partial derivatives $\f_jq$ are $q_j$.  Then a standard
result about distributional derivatives,~\cite[Thm~1.4.2,
p.~10]{Hormander:book63}, implies that the classical partial
derivatives $\f_jq$ of $q$ are equal to $q_j$ in all of $U$.  But a
function with continuous partial derivatives is $C^1$.  Finally
$\f_jq=q_j$ so the derivative is Lipschitz, that is $q$ is of class
$C^{1,1}$.
\end{proof}

\begin{prop}\label{sym-hess}
Let $p\cn \s^{n-1}\to \R$ be a $C^{1,1}$ function.  Then for almost
all $u\in \s^{n-1}$ the second derivative $\con^2p(u)$ exists and is
self-adjoint.  
\end{prop}

\begin{proof}
If $p$ is $C^{1,1}$ the vector field $\con p$ is Lipschitz and thus 
by Rademacher's Theorem $\con^2p(u)$ exists for almost all $u$.  We
have seen, Remark~\ref{rmk:sym}, that if $p$ is of class $C^2$, then
$\con^2p(u)$ is self-adjoint for all $u\in \s^{n-1}$.  
In the case that $p$ is $C^{1,1}$, for each $\e>0$ there is a $C^2$
function $p_\e$ such that if $E_\e:=\{u\in\s^{n-1}: p(u)=p_\e(u), \con
p(u)=\con p_\e(u), \con^2p(u)=\con^2p_\e(u)\}$ then the measure of
$\s^{n-1}\setminus E_\e$ is less than $\e$, \cite[Thm~3.1.15,
p.~227]{Federer:book}.  As $p_\e$ is $C^2$,
$\con^2p(u)=\con^2p_\e(u)$ is self-adjoint for all $u\in E_\e$.
Letting $\e$ go to zero shows that $\con^2p$ is self-adjoint almost
everywhere on $\s^{n-1}$.
\end{proof}

Before applying Proposition~\ref{Gauss-def} to the support function of
a convex set, it is useful to record some symmetry properties of the
operators $\con$ and $\con^2$.  Note that the tangent spaces
$T_u\s^{n-1}$ and $T_{-u}\s^{n-1}$ to $\s^{n-1}$ at antipodal points
$u$ and $-u$ are both just the orthogonal compliment $u^\bot$ to $u$.
Therefore for a function $p$ on $\s^{n-1}$ the vectors $\con p(u)$ and
$\con p(-u)$ are in the same vector space, $u^\bot$, and the linear
maps $\con^2p(u)$ and $\con^2p(-u)$ act on the same vector space
$u^\bot$.  Recall that a function $p\cn \s^{n-1}\to\R$ is even
(respectively odd) iff $p(-u)=p(u)$ (respectively $p(-u)=-p(u)$).
These definitions extend in a obvious way to vector fields or fields
of linear maps on $\s^{n-1}$.  The proof of the following is elementary
and left to the reader.

\begin{lemma}\label{odd}
Let $p\cn \s^{n-1}\to \R$ be a $C^{1,1}$ function.  If $p$ is  even,
then $\con p$ is odd, and $\con^2p$ is even.  If $p$ is odd, then
$\con p$ is even, and $\con^2p$ is odd.  (As $p$ is $C^{1,1}$ the
tensor $\con^2p$ will only be defined almost everywhere.  Saying this
it is even (or odd) means that $\con^2p(u)$ is defined if and
only if $\con^2(-u)$ is defined and at these points
$\con^2p(-u)=\con^2p(u)$ (or $\con^2p(-u)=-\con^2p(u)$).)\qed
\end{lemma}


Recall that if $K$ is a convex body with $C^1$ boundary $\f K$, then
the \emph{Gauss map} is the function $\nu \cn\f K\to \s^{n-1}$ where
$\nu(x)=u$ iff $u$ is the (unique as $\f K$ is $C^1$) outward pointing
unit vector to $K$ at $x$.  If $h$ is the support function of $K$,
then it is not hard to see that $h(\nu(x))=\la x,\nu(x)\ra$,
Therefore, if $\nu$ is injective, so that $\nu^{-1}$ exists, then
$h(u)=\la \nu^{-1}(u),u\ra$, \cite[p.~106]{Schneider:convex}.  More
generally when the support function $h$ is $C^1$ the function
$\phi(u)=h(u)u+\con h(u)$ can still be viewed as the inverse of the
Gauss map:

\begin{prop}\label{spt-form}
Let $K$ be a convex body in $\R^n$ with $C^1$ support function $h$.
Then the map $ \phi(u)=h(u)u+\con h(u) $ maps $\s^{n-1}$ onto $\f K$
with the property that $\phi(x)=u$ if and only if $u$ is an outward
unit normal to $K$ at~$x$.  
\end{prop}

\begin{proof}
We first assume that $\f K$ is $C^\infty$ with positive curvature.
Then the Gauss map $\nu\cn \f K\to \s^{n-1}$ is a diffeomorphism.  Let
$\phi:=\nu^{-1}\cn \s^{n-1}\to \f K$ be the inverse of $\nu$.  Then
$\phi$ is a diffeomorphism and $T_u\s^{n-1}$ and $T_{\phi(u)}\f K$ are
the same (as we are identifying subspaces that differ by a parallel
translation).  Whence $\phi'(u)X\in T_u\s^{n-1}$ for all $X\in
T_u\s^{n-1}$.  By Proposition~\ref{spt-form} this implies there is a
unique smooth real valued function $p$ on $\s^{n-1}$ such that
$\phi(u)=p(u)u+\con p(u)$.  Then $p(u)=\la \phi(u),u\ra$.  But, from
the remarks above, the support function of $K$ is also given by
$h(u)=\la \phi(u),u\ra$ and therefore $p=h$.  So in this case
$\phi(u)=h(u)u+\con h(u)$ is the inverse of the Gauss map and so
$\phi(u)=x$ if and only if $u$ is the outward normal to $K$ at $x$ is
clear.
	
Now assume that $h$ is $C^1$ and set $\phi(u)=h(u)u+\con h(u)$.  Then
$\phi$ is a continuous map from $\s^{n-1}$ to $\R^n$.  There are
convex bodies $\{K_\ell\}_{\ell=1}^\infty$ whose boundaries are smooth
with positive curvature and such that if the support function of
$K_\ell$ is $h_\ell$, then $h_\ell \to h$ in the $C^1$ topology,
\cite[pp. 158--160]{Schneider:convex}.  Therefore if
$\phi_\ell(u):=h_\ell(u)u+\con h_\ell(u)$, then $\phi_\ell\to \phi$
uniformly.  The Hausdorff distance (see \cite[p.~48]{Schneider:convex}
for the definition) between $K$ and $K_\ell$ is given in terms of the
support functions by $d_\text{\rm
Hau}(K,K_\ell)=\|h-h_\ell\|_{L^\infty}$, \cite[1.8.11,
p.~53]{Schneider:convex}, and so $K_\ell \to K$ in the Hausdorff
metric.  Because $K$ and $K_\ell$ are convex this implies $\f
K_\ell\to\f K$ in the Hausdorff metric.  As $\phi_\ell(u)\in \f
K_\ell$ this yields $\phi(u)=\lim_{\ell \to \infty}\phi_\ell(u)\in \f
K$.  Therefore $\phi$ maps $\s^{n-1}$ into $\f K$.  Let $x\in \f K$
and let $u$ be an outward pointing unit normal to $K$ at $x$.  Then
$u$ is an outward pointing normal to $K_\ell$ at $\phi_\ell(u)$.
Therefore the half space $H^-_\ell:=\{y\in \R^n: \la y,u\ra \le
h_\ell(u)\}$ contains $K_\ell$ and its boundary $\f H^-_\ell$ is a
supporting hyperplane to $K_\ell$ at $\phi_\ell(u)$.  Using that
$h_\ell \to h$ uniformly, that $K_\ell\to K$ in the Hausdorff metric,
and that $\phi_\ell(u)\to \phi(u)$ we see that $K$ is contained in
$H^+:=\{y\in \R^n: \la y,u\ra \le h(u)\}$ and that $x\in \f H^+$.
Thus $u$ is an outward pointing unit normal to $K$ at $\phi(u)$.  But,
\cite[Cor.~1.7.3, p.~40]{Schneider:convex}, if the support function is
differentiable, then the body is strictly convex.  Therefore $K$ is
strictly convex and thus a unit vector can be an outward unit normal
to $K$ in at most one point.  So, as $u$ is an outward unit normal to
$K$ at $\phi(u)$ and at $x$, we have $\phi(u)=x$.

Summarizing, if $x\in \f K$ and $u$ is an outward unit normal to $K$
at $x$, then $\phi(u)=x$.  But for any point of $\f K$ there is at
least one unit normal $u$ to $K$ at $x$, so $\phi\cn \s^{n-1}\to \f K$
is surjective.  To finish we need that if $\phi(u)=x$, then $u$ is an
outward pointing unit normal to $K$ at $x$.  The vector $u$ will be an
outward pointing unit normal to $K$ at some point $y\in \f K$.  But
then $\phi(u)=y$.  Thus $x=y$ and $u$ is an outward pointing unit
vector to $K$ at $x$.
\end{proof}

\begin{prop}\label{low-area}
Let $K$ be a compact body with $C^{1,1}$ support function $h$.  Then
$hI+\con^2h$ is positive semi-definite almost everywhere on
$\s^{n-1}$.  If in addition there is a Euclidean ball that slides
freely inside of $\s^{n-1}$, then there is a positive constant $C_1$
such that $\det(hI+\con^2h)\ge C_1$ almost everywhere on $\s^{n-1}$.
\end{prop}

\begin{proof}
Let $\tilde{h}$ be the extension of $h$ to $\R^n$ as a homogeneous
function of degree one (thus $\tilde{h}$ is given by both the
formulas\eq{extend1} and\eq{extend2}).  The function $\tilde{h}$ is
convex, \cite[Thm~1.7.1, p.~38]{Schneider:convex}, and therefore its
Hessian $\f^2\tilde{h}$ is positive semi-definite at all points where
it exists and is self-adjoint.  But then the formula\eq{f2-con2}
relating $\f^2\tilde{h}$ and $\con^2h$ together with
Remark~\ref{rmk:sym} and Proposition~\ref{sym-hess}, shows that
$hI+\con^2h$ is positive semi-definite almost everywhere on
$\s^{n-1}$.

Assume that the Euclidean ball $B_{2r}$ of radius $2r$ slides freely
inside of $K$.  Then there is a convex set $K_1$ such that
$K_1+B_{2r}=K$.  However $K_1$ may not be a convex body.  But
$K_1+B_{2r}=(K_1+B_r)+B_r$ and $K_1+B_r$ is a convex body.  So by
replacing $K_1$ by $K_1+B_r$ we can assume $K_1+B_r=K$ with $K_1$ a
convex body.  Let $h_1$ be the support function of $K_1$.  Then, as
the support function of $B_r$ is the constant $r$, $h_1+r=h$.  This
implies that $h_1$ is also $C^{1,1}$ and therefore $(h_1I+\con^2h_1)$
is positive semi-definite almost everywhere.  But for any positive
semi-definite matrices $A$ and $B$ the inequality $\det(A+B)\ge
\det(A)$ holds.  Therefore
$$
\det(hI+\con^2h)=\det(rI+(h_1I+\con^2h_1))\ge \det(rI)=r^{n-1}=:C_1.
$$
almost everywhere.  
\end{proof}

\begin{lemma}\label{proj-area}
Let $K$ be a convex body in $\R^n$ with $C^{1,1}$ support function
$h$.  Then for any unit vector $a\in \R^n$,
$$
2V_{n-1}(K|a^\bot)=\int_{\s^{n-1}}\det(hI+\con^2 h)|\la a,u\ra|\,dV_{n-1}(u).
$$
\end{lemma}

\begin{proof}
Let $h$ be the support function of $K$ and let $\phi\cn \s^{n-1}\to \f
K$ be $\phi(u)=h(u)u+\con h(u)$.  By Proposition~\ref{spt-form} $\phi$
maps $\s^{n-1}$ onto $\f K$ and, as $h$ is $C^{1,1}$, the map $\phi$
is Lipschitz.  As $\phi$ is Lipschitz it is differentiable almost
everywhere and by Proposition~\ref{Gauss-def} at the points $u$ where
it is differentiable $\phi'(u)=h(u)I+\con^2h(u)$.  Let $f\cn
\s^{n-1}\to K|a^\bot$ be the function $f(u)=\phi(u)|a^\bot$.  This
maps $\s^{n-1}$ onto $K|a^\bot$.  An elementary computation shows that
the Jacobian, $J(f)(u):=\det(f'(u))$, of $f$ is given by
$J(f)(u)=\det\(h(u)I+\con^2h(u)\)\la a,u\ra$.  The area
theorem,~\cite[Thm.~3.2.3, p.~243]{Federer:book}, (note that the
definition of Jacobian used in the area theorem is the absolute value
of the one being used here) implies
\begin{align*}
\int_{K|a^\bot} \# (f^{-1}[y])\,dV_{n-1}(y)&=
\int_{\s^{n-1}}|J(f)(u)|\,dV_{n-1}(u)\\
&=\int_{\s^{n-1}}\det\(h(u)I+\con^2h(u)\)|\la a,u\ra|\,dV_{n-1}(u)
\end{align*}
where $\# (f^{-1}[y])$ is the number of points in the preimage
$f^{-1}[y]:=\{x: f(x)=y\}$.  To complete the proof it is enough to
show $\# (f^{-1}[y])=2$ for almost all $y\in K|a^\bot$.

As $K|a^\bot$ is convex its boundary $\f (K|a^\bot)$ has measure zero.
Therefore we only need consider $y$ in the interior, $\interior\(
K|a^\bot\)$, of $K|a^\bot$.  If $y\in\interior\( K|a^\bot\)$ then
there are exactly two points $x_1,x_2\in \f K$ with $x_j|a^\bot=y$.
Thus $f^{-1}[y]$ is the disjoint union of $\phi^{-1}[x_1]$ and
$\phi^{-1}[x_2]$.  But, \cite[Thm~2.2.4, p.~74]{Schneider:convex}, the
set, $P$, of points $x$ in $\f K$ such that there is more than one
outward unit normal to $K$ at $x$ is a set of measure zero.  So if
$x_1,x_2\notin P$, each of the sets $\phi^{-1}[x_1]$ and
$\phi^{-1}[x_2]$ will have just one elment and therefore
$\#(f^{-1}[y])=2$.  The map $y \mapsto y|a^\bot$ is Lipschitz and
therefore it maps sets of measure zero to sets of measure zero.  Thus
$P|a^\bot$ is a set of measure zero.  Whence for $y\in \interior\(
K|a^\bot\)\setminus P|a^\bot$, and therefore for almost all $y\in
K|a^\bot$, $\#(f^{-1}[y])=2$ which finishes the proof.
\end{proof}

\begin{prop}\label{same:bright}
Let $K_1$ and $K_2$ be convex bodies in $\R^n$ with $C^{1,1}$ support
functions $h_1$ and $h_2$ respectively.  Then there is a constant
$\beta$ such that $V_{n-1}(K_1|a^\bot)=\beta V_{n-1}(K_2|a^\bot)$
for all $a\in \s^{n-1}$ if and only if
$$
\det(h_1I+\con^2h_1)=\beta \det(h_2I+\con^2 h_2)+q,\quad\text{with
$q$ an odd function.}
$$
\end{prop}

\begin{proof}
By Lemma~\ref{proj-area} $V_{n-1}(K_1|a^\bot)=\beta
V_{n-1}(K_2|a^\bot)$ for all $a\in \s^{n-1}$ if and only if
$\int_{\s^{n-1}}q(u) |\la a,u\ra|\,du=0$ for all $a\in \s^{n-1}$ where
$q=\det(h_1I+\con^2h_1)-\beta \det(h_2I+\con^2 h_2)$.  That is, if and
only if $q$ is in the kernel of the cosine transform $(Cf)(a):=
\int_{\s^{n-1}}f(u)|\la a,u\ra|\,du$.  But, \cite[Thm~C,2.4,
p.~381]{Gardner:book}, the kernel of the cosine transform is exactly
the set of odd functions on $\s^{n-1}$.
\end{proof}

\section{Three dimensional Bodies of Constant Width and
Brightness.}\label{sec:3D} 

To prove Theorem~\ref{main:rel} we let $K$ and $K_0$ be convex bodies
in $\R^3$ such that $K_0$ is centrally symmetric about the origin and
that there are constants $\alpha$ and $\beta$ such that $
w_K(u)=\alpha w_{K_0}(u)$ and $(K|u^{\bot})=\beta V_{2}(K_0|y^\bot) $
for all unit vectors $u$.  By rescaling $K$ by a factor of $1/\alpha$
we can assume that $\alpha=1$, that is $K$ and $K_0$ have same width
in all directions.  Then $K_0$ being centrally symmetric about the
origin implies that $K_0$ is the central symmetral $\frac12(K-K)$ of
$K$.  Therefore to prove Theorems~\ref{thm:main} and~\ref{main:rel} it
is enough to prove:

\begin{thm}\label{thm:norm}
Let $K$ be a convex body in $\R^3$ such that its central symmetral
$K_0=\frac12(K-K)$ is a regular gauge and for some constant~$\beta$
\begin{equation}\label{brt}
V_{2}(K|u^{\bot})=\beta V_{2}(K_0|u^\bot)\quad\text{for all $u\in\s^2$}
\end{equation} 
Then $K$ is a translate of $K_0$.
\end{thm}

\begin{lemma}\label{lambda}
If\eq{brt} holds, then $\beta\le 1$ and if $\beta=1$, then $K$ is
a translate of $K_0$.
\end{lemma}

\begin{proof}
Let $u\in \s^2$.  Then $K_0|u^\bot$ is centrally symmetric about the
origin and, viewed as convex bodies in the two dimensional space
$u^\bot$, the sets $K_0|u^\bot$ and $K|u^\bot$ have the same width
function.  Therefore $K_0|u^\bot$ is the central symmetral of
$K|u^\bot$.  By Proposition~\ref{K0-vol} this implies
$V_2(K_0|u^\bot)\ge V_2(K|u^\bot)$ with equality if and only if
$K|u^\bot$ is a translate of $K_0|u^\bot$.  As $V_2(K|u^\bot)=\beta
V_2(K_0|u^\bot)$ this yields that $\beta\le1$.  If $\beta=1$, then for
all $u\in\s^2$ the set $K|u^\bot$ is a translate of $K_0|u^\bot$.
This implies, \cite[Thm~3.1.3, p.~93]{Gardner:book}, that $K$ is a
translate of $K_0$.
\end{proof}

From now on we assume $K$ and $K_0$ satisfy the hypothesis of
Theorem~\ref{thm:norm} and that $h$ and $h_0$ are the support
functions of $K$ and $K_0$ respectively.  By Lemma~\ref{lambda} if
$\beta=1$, Theorem~\ref{thm:norm} holds, so, towards a contradiction,
assume $\beta<1$.

\begin{lemma}\label{beta<1}
If $\beta<1$ then $h$ and $h_0$ are related by $ h=h_0+p$ where $p$ is
an odd function.  The function $p$ satisfies
\begin{enumerate}
\item $p$ is of class $C^{1,1}$,
\item The equality 
\begin{equation}\label{det-p}
\det(pI+\con^2p)=-(1-\beta)\det(h_0I+\con^2h_0)
\end{equation}
holds almost everywhere on $\s^{n-1}$.  Therefore there is a constant
$\delta_0>0$ such that
\begin{equation}\label{p-ineq}
\det(pI+\con^2p)\le -\delta_0
\end{equation}
almost everywhere on $\s^{n-1}$.
\item If $\phi\cn \s^{2}\to \R^3$ is given by $ \phi(u)=p(u)u+\con
p(u) $ then $\phi$ is Lipschitz and $\phi(-u)=\phi(u)$.
\end{enumerate}
\end{lemma}

\begin{proof}
As $K$ and $K_0$ have the same width function,
$h(u)+h(-u)=h_0(u)+h_0(-u)=2h_0(u)$ as $h_0(-u)=h_0(u)$ because $K_0$ is
centrally symmetric about the origin.  Therefore
$$
h(u)=\frac12(h(u)+h(-u))+\frac12(h(u)-h(-u))=h_0(u)+p(u)
$$
where $p(u):=\frac12(h(u)-h(-u))$ is clearly an odd function. 

As $K_0$ is a regular gauge it slides freely inside of some Euclidean
ball and thus by Proposition~\ref{Hug:slide} $h_0$ is $C^{1,1}$.  Then
Corollary~\ref{even-reg} implies $h$ is $C^{1,1}$ and the formula
$p(u)=\frac12(h(u)-h(-u))$ shows that $p$ is also $C^{1,1}$.  

Proposition~\ref{same:bright} implies there is an odd function $q$ on
$\s^2$ such that
\begin{equation}\label{det-h}
\det(hI+\con^2h)=\beta\det(h_0I+\con^2h_0)+q
\end{equation}
holds almost everywhere on $\s^2$.  The equality $h=h_0+p$ implies
\begin{equation}\label{det-h-1}
\det(hI+\con h)=\det\((pI+\con^2p)+(h_0I+\con^2h_0)\).
\end{equation}
For any $2\times 2$ matrix $\trace(A)^2-\trace(A^2)=2\det(A)$, where
$\trace(A)$ is the trace of $A$.  Define $\sigma(A,B)$ on pairs of
$2\times 2$ matrices by
$\sigma(A,B)=\frac12\(\trace(A)\trace(B)-\trace(AB)\)$.  Then
$\sigma(\,,)$ is a symmetric bilinear form and $\sigma(A,A)=\det(A)$.
Whence $\det(A+B)=\det(A)+2\sigma(A,B)+\det(B)$.  Using this
in\eq{det-h-1} gives
\begin{equation}\label{det-h-2}
\det(hI+\con^2h)=\det(pI+\con^2p)+2\sigma(pI+\con^2p,h_0I+\con^2h_0)
+\det(h_0I+\con^2h_0).
\end{equation}
The function $h_0$ is even on $\s^2$ and Lemma~\ref{odd} implies
$\con^2h_0$ is also even.  Therefore $h_0I+\con^2h_0$ is even.
Likewise Lemma~\ref{odd} applied to the odd function $p$ implies
$pI+\con^2p$ is odd.  But $\det(-A)=\det(A)$ for $2\times 2$ matrices,
so the function $\det(pI+\con^2p)$ is even.  The function
$\sigma(pI+\con^2p,h_0I+\con^2h_0)$ is odd as a function of the first
argument and even as a function of the second argument, therefore
$\sigma(pI+\con^2p,h_0I+\con^2h_0)$ is an odd function.  Comparing the
two formulas\eq{det-h} and\eq{det-h-2} for $\det(hI+\con^2h)$ and
equating the even parts gives
$$
\beta\det(h_0I+\con^2h_0)=\det(pI+\con^2p)+\det(h_0I+\con^2h_0).
$$
This implies\eq{det-p}.  By Proposition~\ref{low-area} and the
assumption that $K_0$ slides freely inside of a Euclidean ball there
is a constant $C_1>0$ such that $\det(h_0I+\con^2h_0)\ge C_1$.
Then\eq{det-p} implies\eq{p-ineq} holds with $\delta_0=(1-\beta)C_1$.

That $p$ is $C^{1,1}$ implies $\phi(u)=p(u)u+\con p(u)$ is Lipschitz.
The function $p$ is odd and, by Lemma~\ref{odd}, the vector field $\con
p$ is even.  Therefore $ \phi(-u)=p(-u)(-u)+\con p(-u)=p(u)u+\con
p(u)=\phi(u)$.
\end{proof}

Letting $p$ and $\phi(u)=p(u)u+\con p(u)$ be as in the last lemma, for
any unit vector $a$ let $H_a:=\la \phi(x),a\ra$ be the height function
of $\phi$ in the direction~$a$.  The following, which is trivial when
$h$ is $C^2$ (so that $\phi$ is $C^1$), is the main geometric fact behind
the proof of Theorem~\ref{thm:norm}.

\begin{claim}\label{claim}
If the height function $H_a$ has a local maximum or minimum at $u_0$,
then $u_0=\pm a$.
\end{claim}

\begin{proof}[Proof of Theorem~\ref{thm:norm} assuming the Claim.]
By compactness of $\s^2$ and the continuity of the height function
$H_a$, there are points $u_1,u_2\in \s^2$ such that $H_a(u_1)$ is a
global minimum and $H_a(u_2)$ is a global maximum of $H_a$.  By the
claim $u_1=\pm a$ and $u_2=\pm a$, and therefore $u_1=\pm u_2$.  By
Lemma~\ref{beta<1}, 
$\phi$ is an even function on $\s^2$ and whence
$$
H_a(u_1)=\la \phi(u_1),a\ra=\la \phi(\pm u_2),a\ra=\la
\phi(u_2),a\ra=H_a(u_2). 
$$
As $H_a(u_1)$ and $H_a(u_2)$ are the minimum and maximum of $H_a$ this
implies $H_a(u)$ is constant.  But this is true for any choice of $a$,
so $\phi$ is constant.  Then $\phi'(u)=0$ for all $u\in \s^2$.
However, by Proposition~\ref{Gauss-def}, $\phi'(u)=p(u)I+\con^2 p(u)$
for almost all $u\in\s^2$ and, by Lemma~\ref{beta<1},
$\det(pI+\con^2p)<0$ almost everywhere, which implies $\phi'(u)\ne0$
for almost all $u$.  This contradiction completes the proof.
\end{proof}

We now reduce the claim to an analytic lemma that is proven in the
next section.  Let $e_1, e_2, e_3$ be the standard basis of $\R^3$.
By a rotation we can assume that the height function, $H_a$, has a local
maximum or maximum at $e_3$.  Then to prove the claim we need to show
that $a=\pm e_3$.  We parameterize the open upper hemisphere $\s^2_+$
of $\s^2$ by
\begin{equation}\label{u-def}
u=u(x,y):=\[\begin{matrix}x\\ y\\ \sqrt{1-x^2-y^2}\end{matrix}\]
\end{equation}
where $(x,y)\in \Delta_1:=\{(x,y): x^2+y^2<1\}$.  The function $p$
restricted to $\s^2_+$ can be expressed in terms of the coordinates
$x,y$. Direct calculation shows
$$
\con p= \[\begin{matrix}p_x\\ p_y\\ 0\end{matrix}\]-\left\la \[\begin{matrix}p_x\\ p_y\\ 0\end{matrix}\],u\right\ra u=
\[\begin{matrix}p_x\\ p_y\\ 0\end{matrix}\]-(xp_x+yp_y)
\[\begin{matrix}x\\ y \\ \sqrt{1-x^2-y^2}\end{matrix}\]
$$
Therefore $\phi$ is given by
$$
\phi(x,y)=pu+\con p=\[\begin{matrix}xp\\yp\\p\sqrt{1-x^2-y^2}\end{matrix}\]
+\[\begin{matrix}p_x\\p_y\\0\end{matrix}\]-(xp_x+yp_y)
\[\begin{matrix}x\\y\\\sqrt{1-x^2-y^2}\end{matrix}\]
$$
and
$$
\la \phi,e_3\ra=(p-(xp_x+yp_y))\sqrt{1-x^2-y^2}.
$$
As $p$ is of class $C^{1,1}$, Taylor's theorem implies
$p(x,y)=p(0,0)+xp_x(0,0)+yp_y(0,0)+O(x^2+y^2)$,
$xp_x(x,y)=xp_x(0,0)+O(x^2+y^2)$, and
$yp_y(x,y)=yp_y(0,0)+O(x^2+y^2)$. Therefore
\begin{equation}\label{e3-est}
\la \phi,e_3\ra=p(0,0)+O(x^2+y^2).
\end{equation}

We also consider the projection of $\phi$ onto the first two
coordinates:
$$
\psi(x,y):=
\[\begin{matrix}p_x+x(p-xp_x-yp_y)\\ p_y+y(p-xp_x-yp_y)\end{matrix}\].
$$
This is clearly Lipschitz in a neighborhood of the origin.

\begin{mlemma}\label{lem:z-w}
With $\psi$ as above, there is an open neighborhood $W$ of $\psi(0,0)$
in $\R^2$ and a constant $C_o$ such that for all $w\in W$ there is a
$z\in \Delta_1$ with $\psi(z)=w$ and $C_o^{-1}|z|\le |w-\psi(0,0)|\le
C_o|z|$.
\end{mlemma}

Assuming this we prove Claim~\ref{claim}.  Write the unit vector $a$
defining the height function $H_a$ as $a=\tilde{a}+a_3e_3$ where
$\tilde{a}\in \R^2$ and $a_3\in \R$.  Then for $z=(x,y)\in \Delta_1$ and
using\eq{e3-est}
\begin{align}
H_a(z)&=\la \phi(z),a\ra = \la \psi(z),\tilde{a}\ra+a_3\la
\phi(z),e_3\ra\nonumber\\
&=\la \psi(z),\tilde{a}\ra+a_3p(0,0)+O(|z|^2).\label{Ha}
\end{align}
For real $t$ with $|t|$ small let $w_t=\psi(0,0)+t\tilde{a}$.  By
Lemma~\ref{lem:z-w} there is a $z_t\in \Delta_1$ with $\psi(z_t)=w_t$ and 
$$
|z_t|\le C_o|w_t-\psi(0,0)|=C_o|\tilde{a}||t|.
$$
Thus $|z_t|^2=O(t^2)$.  Using this in\eq{Ha} gives
\begin{align*}
H_a(z_t)&=\la \psi(0,0)+t\tilde{a},\tilde{a}\ra+a_3p(0,0)+ O(|z_t|^2)\\
&=\(\la \psi(0,0), \tilde{a}\ra+a_3p(0,0)\) +t|\tilde{a}|^2+ O(t^2)
\end{align*}
This can only have a local maximum or minimum at $t=0$ if
$\tilde{a}=0$.  As $a$ is a unit vector this implies that $a=\pm e_3$ 
and completes the proof of Claim~\ref{claim}.

\section{Quasiconformal maps and the proof of the main lemma.}\label{sec:quasi}

\subsection{Preliminaries on quasiconformal maps and the Beltrami
equation.}  We recall some basic definitions and facts about
quasiconformal maps.  We identify the complex numbers $\C$ with the
real plane $\R^2$.  Let $U\subseteq \C$ be an open set.  If $f\cn U\to
\C$ write $f=u+iv$.  The function $f$ is in the Sobolev space
$W^{1,2}_\text{\rm Loc}(U)$ iff its distributional first derivatives
are measurable functions that are square integrable on any compact
subset of $U$.  If $f\in W^{1,2}_\text{\rm Loc}(U)$, then the partial
derivatives $f_x=u_x+iv_x$ and $f_y=u_y+iv_y$ exist almost everywhere
on $U$.  If $f\cn U\to V$ is a homeomorphism between open sets $U$ and
$V$ and also $f\in W^{1,2}_\text{\rm Loc}(U)$ a theorem of Gehring and
Lehto~\cite{Gehring-Lehto:diff} (cf.~\cite[Lem.~1,
p.~24]{Ahlfors:quasi}) implies that $f$ is differentiable almost
everywhere (where the derivative, $f'(z)$, is a real linear map
$f'(z)\cn \R^2\to \R^2$).  The \emph{operator norm} of the linear map
$f'(z)$ is $\|f'(z)\|:=\sup_{|v|=1}|f'(z)v|$ and the Jacobian is
$J(f)(z)=\det(f'(z))=u_xv_y-u_yv_x$. For $K\ge 1$ a homeomorphism
$f\cn U\to V$ between two open subsets of $\C$ is
$K$-\emph{quasiconformal} iff $f\in W^{1,2}_\text{\rm Loc}(U)$ and
$$
\|f'(z)\|^2\le KJ(f)(z)
$$
holds almost everywhere in $U$.  There are other equivalent analytic
definitions of $K$-quasiconformality (cf.~\cite[p.~24]{Ahlfors:quasi},
\cite[pp.~6--7]{Krushkal:quasi}, \cite[p.~5]{Iwaniec-Martin:book}).
Introducing the complex derivatives $ {\f}/{\f z}=\frac12\({\f}/{\f
x}-i{\f}/{\f y}\)$ and ${\f}/{\f \ol{z}}=\frac12\({\f}/{\f
x}+i{\f}/{\f y}\) $ an equivalent definition for a homeomorphism $f\cn
U\to V$ which is in $W^{1,2}_\text{\rm Loc}(U)$ to be
$K$-quasiconformal is that
$$
|f_{\ol{z}}|\le \frac{K-1}{K+1}|f_z|
$$
almost everywhere on $U$.  There is a geometric definition of
$K$-quasiconformal (for example see \cite[p.~21]{Ahlfors:quasi}) that
has the advantage that it makes it clear that a homeomorphism $f\cn
U\to V$ between open subsets of $\C$ is $K$-quasiconformal if and only
if its inverse $f^{-1}\cn V\to U$ is $K$-quasiconformal.  The
equivalence of the geometric and analytic definitions of
$K$-quasiconformal was proven by Gehring and Lehto
in~\cite{Gehring-Lehto:diff} (cf.~\cite[Chap.~II]{Ahlfors:quasi}).  A
corollary of the Gehring and Lehto theorem is the following (which can
also be found explicitly in~\cite[Thm~4, p.~9]{Krushkal:quasi}).

\begin{prop}\label{prop:inverse}
If $f\cn U\to V$ is a $K$-quasiconformal map between open subsets of
$\C$, then the inverse $f^{-1}\cn V\to U$ is also $K$-quasiconformal
and satisfies
$$
(f^{-1})_{w}
=\frac{\ol{f_z}}{|f_z|^2-|f_{\ol{z}}|^2},
\qquad
(f^{-1})_{\ol{w}}
=\frac{-f_{\ol{z}}}{|f_z|^2-|f_{\ol{z}}|^2}.
$$
almost everywhere on $V$.\qed
\end{prop}

This implies a result on the Lipschitz invertiblity of certain
homeomorphisms. Let $A\ge 1$, then an open connected subset $V$ of
$\C$ has \emph{$A$-uniformly bounded intrinsic distances} iff any two
points $w_0,w_1\in V$ can be joined by a smooth curve $c$ contained in
$V$ with $\length(c)\le A|w_1-w_0|$.

\begin{prop}\label{lip:inverse}
Let $f\cn U\to V$ be a homeomorphism between open connected subsets of
$\C$ such that the distributional first derivatives of $f$ are bounded
measurable functions and such that the Jacobian satisfies $J(f)\ge
\delta$ almost everywhere for some positive constant $\delta$.  Also
assume $V$ has $A$-uniformly bounded intrinsic distances for some
$A\ge 1$.  Then the inverse $f^{-1}\cn V\to U$ is Lipschitz.
\end{prop}

\begin{lemma}\label{diff-lip}
Let $V$ be an open set in $\C$ with $A$-uniformly bounded intrinsic
distances.  Let $g\cn V\to C$ be a function whose distributional first
derivatives are bounded measurable functions.  Then $g$ is Lipschitz.
\end{lemma}

\begin{proof} We start by constructing the standard smoothing of $g$
by convolution.  Let $\rho$ be a $C^\infty$ non-negative real valued
function on $\C$ with its support contained in the unit disk and with
$\int_\C\rho(s)\,dV_2(s)=1$.  Set $\rho_\e(s):=\e^{-2}\rho(s/\e)$.
Then $\int_\C\rho_\e(s)\,dV_2(s)=1$ and $\rho_\e$ has its support in
the disk of radius $\e$ about the origin.  Let $g_\e(w)=\int_\C
g(w-s)\rho_\e(s)\,dV_2(s)$ be the convolution of $g$ and $\rho_\e$.
Letting $V_\e$ be the set of points in $V$ that are a distance of at
least $\e$ from the boundary $\f V$, $g_\e$ is $C^\infty$ in $V_\e$
and $g_\e \to g$ uniformly on compact subsets of $V$ as $\e\to0$.
Convolution commutes with taking distributional partial
derivatives,~\cite[Thm~1.6.1 p.~14]{Hormander:book63}, and therefore
$$
(g_\e)_x(w)=\int_\C g_x(w-s)\rho_\e(s)\,dV_2(s),\quad
(g_\e)_y(w)=\int_\C g_y
(w-s)\rho_\e(s)\,dV_2(s).
$$
By assumption there is a constant $C_2$ such that $|g_x|,|g_y|\le C_2$
on $V$.  The formulas for $(g_\e)_x$ and $(g_\e)_y$ then show that
$|(g_\e)_x|, |(g_\e)_y|\le C_2$ on $V_\e$.  This implies the operator
norm of $(g_\e)'$ satisfies $\|(g_\e)'(w)\|\le 2C_2$ on $V_\e$.  Let
$w_0,w_1$ be in $V$.  Then there is a smooth curve $c\cn [0,1]\to V$ with
$c(0)=w_0$ and $c(1)=w_1$ and with $\length(c)\le A|w_1-w_2|$.  
For any $\e$ less than the distance of $c$ from the boundary $\f V$ we
have
\begin{align*}
|g_\e(w_1)-g_\e(w_0)|&=
	\left|\int_0^1\frac{d}{dt}g_\e(c(t))\,dt\right|
\le \int_0^1\|(g_\e)'(c(t)\|\,|c'(t)|\,dt\\
&\le 2C_2\length(c)
\le 2 C_2A |w_1-w_0|.
\end{align*}
Taking the limit as $\e\to0$ gives $|g(w_1)-g(w_0)|\le
2C_2A|w_1-w_0|$ and thus  $g$ is Lipschitz as required.
\end{proof}

\begin{proof}[Proof of Proposition~\ref{lip:inverse}]
Let $z\in U$ and let $r>0$ be small enough that the disk $B(z,r)$ is
contained in $U$.  The restriction of $f$ to $B(z,r)$ will still have
bounded distributional first derivatives and $B(z,r)$ is convex,
therefore Lemma~\ref{diff-lip} implies that $f\big|_{B(z,r)}$ is
Lipschitz.  This shows that $f$ is locally Lipschitz on $U$.  Thus by
Rademacher's Theorem its derivative $f'(z)$ exists almost everywhere
on $U$.  For a locally Lipschitz function the ordinary first partial
derivatives are the some as the distributional first partial
derivatives, whence the assumption about $f$ having bounded first
distributional derivatives implies there is a constant $C_3$ such that
$\|f'(z)\|\le C_3$ almost everywhere on $U$.  But then $\|f'(z)\|^2\le
(C_3^2/\delta)\delta\le (C_3^2/\delta) J(f)$ almost everywhere.
Therefore $f$ is $K$-quasiconformal with $K=(C_3^2/\delta)$.
Calculation shows that the Jacobian is given by
$J(f)=|f_z|^2-|f_{\ol{z}}|^2$ and that $|f_z|, |f_{\ol z}|\le
\|f'(z)\|\le C_3$.  Combining this with Proposition~\ref{prop:inverse}
yields that the distributional derivatives $(f^{-1})_{w}$ and
$(f^{-1})_{\ol{w}}$ are functions with
$$
|(f^{-1})_{w}|
\le \frac{|\ol{f_z}|}{\delta}\le \frac{C_3}\delta,\qquad
|(f^{-1})_{\ol{w}}|
\le \frac{|f_{\ol{z}}|}{\delta} \le \frac{C_3}\delta.
$$
Therefore the distributional first derivatives of $f^{-1}$ are
bounded on $V$ and $V$ has $A$-uniformly bounded intrinsic distances.
Whence Lemma~\ref{diff-lip} implies that $f^{-1}$ is Lipschitz.
\end{proof}

Some basic facts about solutions to the Beltrami equation will also be
needed.  Let $U$ be a open subset of $\C$ and $\mu\cn U\to \C$ a
measurable function with $\|\mu\|_{L^\infty}<1$.  Then the
\emph{Beltermi equation} determined by $\mu$ is
$$
f_{\ol{z}}=\mu f_z.
$$
When $\mu\equiv 0$ this is just the Cauchy-Riemann equations.  The
following summarizes the basic facts about existence and uniqueness of
solutions to Beltrami equations and is a combination of a special case
of a basic existence result of C.~B.~Morrey~\cite{Morrey:quasi}
and a factorization theorem of Stoilow.  A good source for these
results is the book \cite{Krushkal:quasi} where \cite[Thm.~2,
p.~8]{Krushkal:quasi} and \cite[Thm.~3, pp.~8--9]{Krushkal:quasi} can
be combined to give:

\begin{thm}\label{thm:Bel-soln}
Let $U$ be a bounded simply connected open subset of $\C$ and $\mu\cn
U\to \C$ a measurable function with $\|\mu\|_{L^\infty}<1$.  Let
$z_0\in U$.  Then there is a quasiconformal map $q\cn U\to \Delta_r$
that satisfies $q_{\ol{z}}=\mu q_z$ and $q(z_0)=0$.  Moreover, if $f\in
W^{1,2}_\text{\rm Loc}(U)$ and satisfies $f_{\ol z}=\mu f_z$ in the
distributional sense in $U$, then $f(z)= \Phi(q(z))$ for a unique
holomorphic function $\Phi$.\qed
\end{thm}

\begin{prop}\label{lem:main}
Let $U$ be an open disk centered at the origin in $\R^2=\C$ and
$f=u+iv$ a Lipschitz function defined on $U$ with $f(0)=0$, and such
that there is a constant $\delta>0$ with
$J(f)=u_xv_y-u_yv_x\ge\delta>0$ almost everywhere.  Then there is an
$r>0$ and a constant $C_o>0$ such that for any $w\in\C$ with $|w|\le
r$ there is a $z\in U$ with $f(z)=w$ and $C_o^{-1}|z|\le |w|\le
C_o|z|$.
\end{prop}

\begin{proof}
By assumption the Jacobian satisfies $J(f)=|f_z|^2-|f_{\ol{z}}|^2\ge
\delta$ almost everywhere.  This implies  $|f_z|^2\ge
\delta+|f_{\ol{z}}|^2\ge \delta$ and thus $|f_z|\ge\sqrt{\delta}>0$
almost everywhere.  Whence the complex valued function
$$
\mu(z)=\frac{f_{\ol{z}}}{f_z}
$$
is defined almost everywhere on $U$.  Also $|f_z|^2-|f_{\ol{z}}|^2\ge
\delta$ implies
$$
|\mu(z)|^2=\left|\frac{f_{\ol{z}}}{f_z}\right|^2\le
1-\frac{\delta}{|f_z|^2}.
$$
But $f$ is Lipschitz so there is a constant $C_4$ with $|f_z|^2\le
C_4$ almost everywhere in $U$.  Using this in the last inequality gives
$$
|\mu(z)|^2\le 1-\frac{\delta}{C_4}:=C_5^2<1.
$$
Thus $f$ satisfies the Beltrami equation $f_{\ol z}=\mu(z)f_z$.
where $\|\mu\|_{L^\infty}\le C_5<1$.

By Theorem~\ref{thm:Bel-soln} there is a homeomorphism $q\cn U\to U$
with $q(0)=0$ and $q\in W^{1,2}_\text{\rm Loc}(U)$ that satisfies $
q_{\ol{z}}=\mu(z)q_z $ and a holomorphic function $\Phi$ defined on
$U$ such that $ f(z)=\Phi(q(z))$.  As $q(0)=0$ and $f(0)=0$ the the
holomorphic function $\Phi$ will have a zero at $0$.  Assume this zero
is of order $k\ge1$.  Then standard results,
\cite[p.~133]{Ahlfors:book3}, about holomorphic maps imply there is a
holomorphic mapping $\Psi$ with $\Psi(0)=0$, which is conformal near
$0$, and such that $\Phi(w)=\Psi(w)^k$.  Then in a neighborhood of $0$
the map $h:=\Psi\circ q$ is a homeomorphism and in this neighborhood
$f(z)=\Psi(q(z))^k=h(z)^k$.  It follows that there is a small positive
real number $r$ such that if $\Delta_r:=\{w: |w|<r\}$ is the disk of
radius $r$ and $\Delta_r^*:=\{w: 0<|w|<r\}$ the pictured disk, and
$U_r$ is the connected component of $f^{-1}[\Delta_r]$ containing $0$,
and $U_r^*=U_r\setminus \{0\}$, then $f\big|_{U_r^*}\cn U_r^*\to
\Delta_r^*$ is exactly $k$ to~$1$, and is in fact a $k$-fold covering
map. (That is each $w\in \Delta_r^*$ has a neighborhood $N$ that is
\emph{evenly covered} in the sense that $f\big|_{U_r^*}^{-1}[N]$ is a
disjoint union of sets $M_1\cd M_k$ such that $f\big|_{U_r^*}$
restricted to each $M_j$ is a homeomorphism of $M_j$ with $N$.) Let
$f_0:= f\big|_{U_r^*}$.  The fundamental groups of $\Delta_r^*$ and
$U_r^*$ are both isomorphic to the additive group of integers $\Z$ and
the image $f_{0*}\[\pi_1(U_r^*)\]$ in $\pi_1(\Delta_r^*)$ is $k\Z$,
the unique subgroup of index $k$ in $\pi_1(\Delta_r^*)$.  Define a map
$\varpi \cn\Delta_r^*\to \Delta_r^*$ by
$$
\varpi(\rho e^{i\theta})=\rho e^{i k\theta}.
$$
This is also a $k$-fold covering map and thus
$\varpi_*\[\pi_1(\Delta_r^*)\]$ also has index $k$ in
$\pi_1(\Delta_r^*)$.  Whence
$\varpi_*\[\pi_1(\Delta_r^*)\]=f_{0*}\[\pi_1(U_r^*)\]$.  Therefore,
\cite[Thm~5.1, p.~156]{Massey:top} or \cite[Thm~5,
p.~76]{Spanier:book}, there is a continuous lifting $\hat{f}_0\cn
U_r^*\to \Delta_r^*$ such that

\centerline{
\begin{picture}(80,60)(0,-5)
\put(0,0){$U_r^*$}
\put(49,0){$\Delta_r^*$}
\put(49,45){$\Delta_r^*$}
\put(15,2.5){\vector(1,0){30}}
\put(15,12){\vector(1,1){30}}
\put(53,42){\vector(0,-1){30}}
\put(28,37){$\hat{f}_0$}
\put(28,8){$f_0$}
\put(56,25){$\varpi$}
\end{picture}
}

%

\noindent
commutes.  Then $f_0=\varpi\circ \hat{f}_0$ and, \cite[Lem.~6.7,
p.~160]{Massey:top} or \cite[Lem.~1, p.~79]{Spanier:book},
$\hat{f}_0$ is also a covering map (which can also easily be checked
from the definitions).  As $f_0=\varpi\circ \hat{f}_0$ and both the
maps $f_0$ and $\varpi$ are $k$~to~$1$ this forces $\hat{f}_0$ to be
$1$~to~$1$.  But a $1$~to~$1$ covering map is a homeomorphism and
thus $\hat{f}_0\cn U_r^*\to \Delta_r^*$ is a homeomorphism.  

In polar coordinates $(\rho,\theta)$ on $\Delta_r^*$ the standard flat
Riemannian metric is given by $g_0:=d\rho^2+\rho^2\,d\theta^2$.  The
pull back of this metric by $\varpi$ is
$\varpi^*g_0=d\rho^2+k^2\rho^2\,d\theta^2$.  Therefore $g_0\le
\varpi^*g_0 \le k^2g_0$.  This shows for any vector $X$ and any point
$z\in \Delta_r^*$ that $|X|\le |\varpi'(z)X|\le k|X|$.  Thus the
operator norms of the linear maps $\varpi'(z)$ and $\varpi'(z)^{-1}$
satisfy
\begin{equation}\label{varpi-ineq}
\|\varpi'(z)\|\le k,\qquad \|\varpi'(z)^{-1}\|\le 1.
\end{equation} 

The map $\varpi\cn \Delta_r^*\to \Delta_r^*$ is $C^\infty$
and\eq{varpi-ineq}, together with the inverse function theorem, shows
that each point $w\in \Delta_r^*$ has a neighborhood $N$ such that
$\varpi\big|_{N} $ is injective, $\varpi[N]$ is a open subset of
$\Delta_r^*$, and $\varpi\big|_{N}$ is a diffeomorphism of $N$ with
$\varpi[N]$.  Let $z_0\in U_r^*$ and let $N$ be such a neighborhood of
$w_0=\hat{f}_0(z_0)$.  The point $z_0$ will have a neighborhood $V$
such that $f_0(z)\in \varpi[N]$ for all $z\in V$.  Thus
$f_0=\varpi\circ \hat{f}_0$ implies
$\hat{f}_0\big|_V=\varpi\big|_{N}^{-1}\circ f_0\big|_V$.  The function
$\varpi\big|_{N}^{-1}$ is $C^\infty$ and $f_0$ is Lipschitz, thus
$\hat{f}_0$ is Lipschitz near~$z_0$.  Therefore $\hat{f}'_0$ exists
almost everywhere in $V$ and for $z\in V$ where $\hat{f}_0'(z)$ exists
use\eq{varpi-ineq} to get
$$
\|\hat{f}_0'(z)\|\le
\left\|\(\varpi^{-1}\big|_{N}^{-1}\)'(f_0(z))\right\|\|f_0'(z)\|
\le \|f_0'(z)\|.
$$
But this holds in a neighborhood of an arbitrary point $z_0$ of
$U_r^*$ and whence $\|\hat{f}_0'(z)\|\le \|f_0'(z)\|$ almost
everywhere on $U_r^*$.  As $f_0$ is Lipschitz there is a constant $C_6$
such that $\|\hat{f}_0'(z)\|\le \|f_0'(z)\|\le C_6$ almost everywhere
on $U_r^*$.  This shows that the distributional first derivatives of
$\hat{f}_0$ are bounded measurable functions.

It is easy to compute that $J(\varpi)=k$.  By assumption, $J(f)\ge
\delta$ and $f_0$ is a restriction of $f$ whence
$$
\delta\le J(f_0)=J(\varpi\circ \hat{f}_0)=J(\varpi)J(\hat{f}_0)
=kJ(\hat{f}_0).
$$
Therefore $J(\hat{f}_0)\ge \delta/k$.  The set $\Delta_r^*$ has
$A$-uniformly bounded intrinsic distances for all $A>1$.  Thus
$\hat{f}_0\cn U_r^*\to \Delta_r^*$ satisfies the conditions of
Proposition~\ref{lip:inverse}.  Whence $\hat{f}_0^{-1}\cn
\Delta_r^*\to U_r^*$ is Lipschitz.

As $\hat{f}_0^{-1}$ is Lipschitz there is a constant $C_7$ such that
for all $w,w_0\in \Delta_r^*$ the inequality
$|\hat{f}_0^{-1}(w)-\hat{f}_0^{-1}(w_0)|\le C_7|w-w_0|$
holds. Therefore for any $z,z_0\in U_r^*$ 
\begin{equation}\label{z-z0}
|z-z_0|=|\hat{f}_0^{-1}(\hat{f}_0(z))-\hat{f}_0^{-1}(\hat{f}_0(z_0))|
\le C_7|\hat{f}_0(z)-\hat{f}_0(z_0)|.
\end{equation}
From the definition of $\varpi$ it clear that $|\varpi(w)|=w$ for all
$w\in \Delta_r^*$.  Thus
$|\hat{f}_0(z_0)|=|\varpi(f_0(z_0))|=|f_0(z_0)|=|f(z_0)|$.  But
$f(0)=0$ and $f$ is continuous and whence $\lim_{z_0\to 0}f(z_0)=0$.
Therefore $\lim_{z_0\to 0}\hat{f}(z_0)=0$ and thus taking the limit as
$z_0\to 0$ in\eq{z-z0} yields
$$
|\hat{f}_0(z)|\ge \frac{1}{C_7}|z|
$$
for all $z\in U_r^*$.  

We now complete the proof of Proposition~\ref{lem:main}.  Let $w\in
\Delta_r^*$.  Then there is a $z\in U_r^*$ with $f(z)=w$.  By the
definition of $f_0$ as the restriction of $f$ we have
$w=f_0(z)=\varpi(\hat{f}_0(z))$.  Again using that
$|\varpi(\xi)|=|\xi|$ we have
$$
|w|=|\varpi(\hat{f}_0(z))|=|\hat{f}_0(z)|\ge \frac{1}{C_7}|z|.
$$
Also, as $f$ is Lipschitz and $f(0)=0$, there is a constant $C_8$ with
$|w|=|f(z)|\le C_8|z|$.  
Letting $C_o=\max\{C_7,C_8\}$ completes the proof.
\end{proof}

\subsection{Proof of the Main Lemma}

We use the notation of the Section~\ref{sec:3D}.  In particular
$\phi(u)=p(u)u+\con p(u)$, $\psi$ is the projection of $\phi$ onto the
first two coordinates and $u=u(x,y)$ is given by\eq{u-def}. 

\begin{lemma}\label{lem:Jac}
There is an open disk $U$ centered at the origin so that for some
constant $\delta>0$ the Jacobian of $\psi$ satisfies $
J(\psi):=\det(\psi')\le -\delta $ almost everywhere in $U$.
\end{lemma}

\begin{proof}
For $(x,y)$ in the unit disk the tangent plane to $\s^2$ at $u(x,y)$
is $u(x,y)^\bot$ and the orientation of this tangent plane is so that
the projection onto the $(x,y)$ plane is orientation preserving.
(This is because $u(x,y)$ is in the upper hemisphere of $\s^2$.)  By
Proposition~\ref{Gauss-def} $\phi'(z)=p(z)I+\con^2p(z)$ almost
everywhere and by Lemma~\ref{beta<1}
$$
J(\phi)=\det(p(z)I+\con^2p(z))\le-\delta_0
$$
for almost all $z$ in the unit disk and for some $\delta_0>0$.  The
projection $\pi$ of the tangent plane $T(\s^2)_u=u^\bot$ onto $\R^2$
has Jacobian $J(\pi)=\la u,e_3\ra$.  As $\psi=\pi\circ \phi$
$$
J(\psi)=J(\pi)J(\phi)=\la u,e_3\ra J(\phi)\le -\la u,e_3\ra\delta.
$$
But $\la u(x,y),e_3\ra=\sqrt{1-x^2-y^2}$ so if $U=\{(x,y):
x^2+y^2<\sqrt{3}/2\}$, then $J(\pi)>1/2$.  Thus on $U$
$J(\psi)<-\delta$ where $\delta=\frac12\delta_0$.
\end{proof}

Returning to the proof of the Main Lemma, let $U$ be as in the last
lemma and let $f\cn U\to \C$ be given by
\begin{equation}\label{f-psi}
f(z)=\ol{\psi(z)}-\ol{\psi(0)}.
\end{equation}
Complex conjugation is an orientation reversing isometry and $\psi$ is
Lipschitz, thus $f$ is also Lipschitz.  The Jacobian of $f$ is $
J(f)=-J(\psi)\ge \delta$.  And clearly $f(0)=0$.  Note that as $f$ and
$\psi$ are related by\eq{f-psi}, then $\psi(z)=w$ if and only if
$f(z)=\ol{w}-\ol{\psi(0)}$.  Therefore the Main Lemma~\ref{lem:z-w}
follows from Proposition~\ref{lem:main}.  This completes the proof of
Theorem~\ref{thm:norm}.

\section*{Acknowledgments}  I am indebted to Daniel Hug for supplying
the statement and a proof of Proposition~\ref{Hug:slide} which greatly
simplified my initial proof of the $C^{1,1}$ regularity of the support
function of a set of constant width.  A remark of Marek Kossowski
lead me to realize the covering space argument in the proof of the
Proposition~\ref{lem:main} was required.  I also had several useful
conversations with Mohammad Ghomi on topics related to this paper.  
%

\providecommand{\bysame}{\leavevmode\hbox to3em{\hrulefill}\thinspace}

\end{document}